\documentclass[11pt]{article}
\usepackage[a4paper]{anysize}\marginsize{3.5cm}{3.5cm}{1.3cm}{2cm}
\pdfpagewidth=\paperwidth \pdfpageheight=\paperheight
\usepackage{amsfonts,amssymb,amsthm,amsmath,eucal}
\usepackage{pgf}
\usepackage{bbm}
\usepackage{tikz} 
\usepackage{mathrsfs}
\usetikzlibrary{arrows}
\usepackage{color}
\usepackage{float}
\usepackage{subfigure}
\usepackage{caption}

\definecolor{qqqqff}{rgb}{0.,0.,0.}

\theoremstyle{plain}
\newtheorem{thm}{Theorem}[section]
\newtheorem{theorem}[thm]{Theorem}

\newtheorem{conjecture}[thm]{Conjecture}

\theoremstyle{definition}
\newtheorem{definition}[thm]{Definition}
\newtheorem{remark}[thm]{Remark}

\newtheorem{problem}[thm]{Problem}

\newtheorem{thevarthm}[thm]{\varthmname}

\newenvironment{varthm*}[1]{\trivlist\item[]{\bf #1.}\it}{\endtrivlist}

\renewcommand\geq{\geqslant}

\renewcommand\leq{\leqslant}

\newcommand\be{\begin{eqnarray*}}
\newcommand\ee{\end{eqnarray*}}

\newcommand\C{\mathbb C}

\renewcommand\P{\mathbb P}

\newcommand\Ha{\mathcal H}

\newcommand\newop[2]{\def#1{\mathop{\rm #2}\nolimits}}
\newop\log{log}
\newop\ord{ord}
\newop\Gal{Gal}
\newop\SL{SL}
\newop\Bl{Bl}
\newop\mult{mult}
\newop\mass{mass}
\newop\div{div}
\newop\codim{codim}
\newop\sing{sing}
\newop\Zeroes{Zeroes}

\newcommand\alphahat{\widehat{\alpha}}

\def\keywordname{{\bfseries Keywords}}%
\def\keywords#1{\par\addvspace\medskipamount{\rightskip=0pt plus1cm
\def\and{\ifhmode\unskip\nobreak\fi\ $\cdot$
}\noindent\keywordname\enspace\ignorespaces#1\par}}
\def\subclassname{{\bfseries Mathematics Subject Classification
(2010)}\enspace}
\def\subclass#1{\par\addvspace\medskipamount{\rightskip=0pt plus1cm
\def\and{\ifhmode\unskip\nobreak\fi\ $\cdot$
}\noindent\subclassname\ignorespaces#1\par}}

\newcommand\beginproof[1]{\trivlist\item[\hskip\labelsep{\em #1.}]}
\newcommand\proofof[1]{\beginproof{Proof of #1}}

\def\endproof{\hspace*{\fill}\endproofsymbol\endtrivlist}

\def\endproofsymbol{\frame{\rule[0pt]{0pt}{6pt}\rule[0pt]{6pt}{0pt}}}

\begin{document}

\author{\L .~Farnik, J.~Gwo\'zdziewicz, B.~Hejmej, M.~Lampa-Baczy\'nska, G.~Malara, J.~Szpond}
\title{Initial sequences and Waldschmidt constants of planar point configurations}
\date{\today}
\maketitle
\thispagestyle{empty}

\begin{abstract}
   The purpose of this work is to extend the classification of planar
   point configurations with low Waldschmidt constants initiated in
   \cite{DST16} and continued in \cite{MosHag16} for all values less than $5/2$.
   As a consequence we prove a conjecture of Dumnicki, Szemberg
   and Tutaj-Gasi\'nska concerning initial sequences with low first differences.
\keywords{symbolic power, point configurations, Waldschmidt constants, Chudnovsky conjecture}
\subclass{14C20 \and 14J26 \and 14N20 \and 13A15 \and 13F20}
\end{abstract}


\section{Introduction}
   In recent years it has become evident that adopting an asymptotic perspective in algebraic geometry
   and in commutative algebra
   often leads to much more regular and clear patterns than obtained when considering isolated phenomena. Examples
   for the fruitfulness of this approach range from questions about base loci of linear series \cite{ELMNP},
   through the growth of higher cohomology groups \cite{dFKL07} to asymptotic syzygies and Betti numbers \cite{EinLaz12}.
   In the present note we are primarily interested in Waldschmidt constants.
   These invariants were introduced by Waldschmidt in \cite{Wal77}
   and the term ``Waldschmidt constant'' was coined some 30 years latter by Dumnicki and Harbourne.
   The interest in Waldschmidt constants stems partly from their applications in describing the
   geometry of the effective cone of algebraic varieties. They are also quite interesting
   invariants studied in their own right, see e.g. \cite{BocFra16}, \cite{BCGHJNSvTT16}, \cite{DHST14}, \cite{GHVT13}. 

   To begin with, recall that for a homogeneous ideal $I$ in the ring of polynomials
   $\C[x_0,\ldots,x_n]$ its \emph{initial degree} $\alpha(I)$ is defined as the
   least integer $t$ such that the homogeneous part of degree $t$ is non-zero.
   If $I$ is a radical ideal in $\C[x_0,\ldots,x_n]$ and $m$ is a positive integer,
   then the $m$-th symbolic
   power $I^{(m)}$ of $I$ equals the $m$-th differential power of $I$ due to the celebrated
   Nagata-Zariski Theorem (see \cite[Corollary 2.9]{SidSul09}). Thus $I^{(m)}$
   is the homogeneous ideal
   consisting of all polynomials vanishing to order at least $m$ along the
   subvariety $\Zeroes(I)\subset \P^n$. Combining these two notions we introduce
   the following object.
\begin{definition}[Initial sequence]
   Let $I$ be a radical ideal in $\C[x_0,\ldots,x_n]$. The \emph{initial sequence}
   of $I$ is the sequence
   $$\alpha(I),\;\alpha(I^{(2)}),\;\alpha(I^{(3)}),\;\ldots.$$
\end{definition}
   It is easy to see that the initial sequence is strictly growing and subadditive, i.e.
   the inequality
   $$\alpha(I^{(m+n)})\leq \alpha(I^{(m)})+\alpha(I^{(n)})$$
   holds for all $m,n\geq 1$.
   Fekete's Lemma (see \cite{Fek23}) implies then the existence of the following limit.
\begin{definition}[Waldschmidt constant]
   The \emph{Waldschmidt constant} of a homogeneous ideal $I\subset\C[x_0,\ldots,x_n]$
   is the real number
   $$\alphahat(I)=\lim\limits_{m\to\infty}\frac{\alpha(I^{(m)})}{m}.$$
\end{definition}
   It follows also from \cite{Fek23} that $\alphahat(I)$ can be alternatively computed as the infimum
   $$\alphahat(I)=\inf \frac{\alpha(I^{(m)})}{m}.$$

   Using the classification types of \cite{GHM09} for $n \leq 8$ points in $\P^2$,
   Harbourne determined the Waldschmidt constant $\alphahat(I_Z)$ for a reduced scheme Z of each type.
   Harbourne's unpublished list is available online at \cite{Har}. In the present work we
   extend this classification a little bit.
   Our first main result shows that the
   geometry prohibits too small Waldschmidt constants of more than
   $9$ points.
\begin{varthm*}{Theorem A}
   Let $Z$ be a set of $s\geq 9$ points in $\P^2$ and let $I$
   be its radical ideal. If
   $$\alphahat(I)< \frac52,$$
   then $Z$
   \begin{itemize}
   \item  is contained in a conic and consequently $\alphahat(I)\leq 2$, or;
   \item  is a set of nine points arranged as indicated in Figure \ref{fig: 9points} and  $\alphahat(I) =\frac{17}{7}$.
   \end{itemize}
\end{varthm*}
   This result in a more detailed form is proved below as Theorem \ref{thm: conf of points}.
   We apply Theorem A to complete the works of Dumnicki,
   Szemberg and Tutaj-Gasi\'nska \cite{DST13} and \cite{DST16}
   concerning properties of initial
   sequences of points in $\P^2$. Our second main result is the proof of Conjecture 4.6
   from \cite{DST13}.
\begin{varthm*}{Theorem B}
   Let $Z$ be a finite set of points in $\P^2$ and let $I$
   be its radical ideal. Assume that for some integer $m$ we have
   $$\alpha(I^{(m+k)})=\alpha(I^{(m)})+2k$$
   for $k=1,2,3,4$. Then $Z$ is contained in a conic.
\end{varthm*}
\begin{remark}
   It is natural to study
   the question if a certain number of differences
   in the initial sequence of $Z$ equal to some integer $d$
   implies that $Z$ is contained in a curve of degree $d$
   for arbitrary $d$. It has been
   proved to be the case for $d=1$ in \cite[Corollary 3.5]{DST13}.
   On the other hand
   \cite[Example 4.14]{DST13} and
   \cite[Proposition 4.2]{DST16} show that for
   $d\geq 3$ the analogous property does not hold.
   Our result completes thus the picture for $d=2$.
   It would be interesting to understand better the role
   of the rationality of lines and conics on $\P^2$ with regard
   to this property.
\end{remark}
   Throughout the paper we work over the field $\C$ of complex numbers.
\section{Initial sequences and their properties}
   Waldschmidt constants of radical
   ideals of points in complex projective spaces
   were studied by Chudnovsky in \cite{Chu81}. In particular he observed that
   the initial degree $\alpha(I)$ of an ideal $I$ (the first term in the initial sequence)
   and its asymptotic cousin $\alphahat(I)$ are surprisingly closely related and asked if
   the following statement is true.
\begin{conjecture}[Chudnovsky]\label{conj:Chudnovsky}
   Let $I$ be a radical ideal of a finite set of points in $\P^n$.
   Then the inequality
   $$\frac{\alpha(I)+n-1}{n}\leq \frac{\alpha(I^{(m)})}{m}$$
   holds for all $m\geq 1$. In particular
   \begin{equation*}\label{eq:inequality alpha alphahat}
      \frac{\alpha(I)+n-1}{n}\leq \alphahat(I).
   \end{equation*}
\end{conjecture}
   This conjecture was proved for $n=2$ by Chudnovsky in \cite{Chu81}
   using some variants
   of Schwarz Lemma and Bombieri results on plurisubharmonic functions.
   Some special cases have been investigated by Demailly in \cite{Dem82}.
   Esnault and Viehweg using methods of complex projective geometry have
   proved the following useful result, see \cite[In\'egalit\'e A]{EsnVie83}.
\begin{theorem}[Esnault -- Viehweg]\label{thm:EV}
   Let $I$ be a radical ideal of a finite set of points in $\P^n$
   with $n\geq 2$. Let $k\leq m$ be two integers. Then
    \begin{equation}\label{inq: wc}
    \frac{\alpha(I^{(k)})+1}{k+n-1}\leq \frac{\alpha(I^{(m)})}{m},
   \end{equation}
   in particular
   $$\frac{\alpha(I^{(k)})+1}{k+n-1}\leq \alphahat(I).$$
\end{theorem}
   In the aforementioned paper \cite{Chu81} Chudnovsky announced
   (in the text right after Theorem 8) the feasibility of enumerating
   all Waldschmidt constants of all configurations of up to $9$ points
   in $\P^2$. In Appendix 1 of that paper, he also listed some
   point configurations and corresponding constants.
   The complete list of all configurations of up to $8$ points
   was found by Geramita, Harbourne and Migliore in
   \cite{GHM09}. They computed all possible Hilbert
   functions of fat point subschemes supported on up to $8$ points in $\P^2$.
   We denote by $\Ha(s,n)$ the configuration number $n$ for $s$ points
   from Tables 1, 3 and 5 in \cite{GHM09}.

   For $9$ or more points there seem to be infinitely many possible
   values of Waldschmidt constants so that a complete classification
   is not possible. Here we show some restrictions on possible values
   of Waldschmidt constants. The following Theorem implies immediately
   Theorem A. Note that configurations of points with $\alphahat(I)<9/4$
   have been classified in \cite[Main Theorem]{DST16}.
\begin{theorem}[Configurations of points with $\alphahat(I)< 5/2$]\label{thm: conf of points}
   Let $I$ be the radical ideal of a set $Z$ of $s$ points in $\P^2$.
   If
   $$\alphahat(I)<\frac52,$$
   then one of the following cases holds:
\begin{itemize}
    \item[a)] $\alphahat(I)=1$ and $Z$ is contained in a line,
    \item[b)] $\alphahat(I)=2-\frac{1}{s-1}$ and exactly $s-1$ points are collinear,
    \item[c)] $\alphahat(I)=2$ and
            \begin{itemize}
            \item $Z$ is contained in a conic, or
            \item $Z$ is the configuration  $\Ha(6,10)$ of six points indicated in Figure \ref{fig: h(6,10)},
            \end{itemize}
                \begin{figure}[H]
                \centering
                \begin{tikzpicture}[line cap=round,line join=round,>=triangle 45,x=1.0cm,y=1.0cm,scale=0.7]
                \clip(2.1,-0.42) rectangle (6.34,4.36);
                \draw [domain=2.1:6.34] plot(\x,{(-4.9576--1.84*\x)/0.14});
                \draw [domain=2.1:6.34] plot(\x,{(--12.4536-2.24*\x)/2.46});
                \draw [domain=2.1:6.34] plot(\x,{(-2.656--0.4*\x)/-2.6});
                \draw [domain=2.1:6.34] plot(\x,{(-5.28153921228557--1.2007559730703197*\x)/-0.42815242376738905});
                \begin{scriptsize}
                \draw [fill=black] (2.74,0.6) circle (2.5pt);
                \draw [fill=black] (2.88,2.44) circle (2.5pt);
                \draw [fill=black] (5.34,0.2) circle (2.5pt);
                \draw [fill=black] (4.268414155426345,0.3648593607036393) circle (2.5pt);
                \draw [fill=black] (3.8402617316589556,1.565615333773959) circle (2.5pt);
                \draw [fill=black] (2.9940423636496676,3.9388424936813435) circle (2.5pt);
                \end{scriptsize}
                \end{tikzpicture}
                \caption{$ $  $\Ha(6,10)$, $\alphahat(I)=2$}\label{fig: h(6,10)}
                \end{figure}
    \item[d)] $\alphahat(I)=9/4$ and $Z$ is the configuration $\Ha(6,9)$ of six points
          indicated in Figure \ref{fig: 94 s=6},
          \begin{figure}[H]
            \centering
               \begin{tikzpicture}[line cap=round,line join=round,x=1.0cm,y=1.0cm,scale=0.7]
               \clip(-4,-0.5) rectangle (4.5,5);
               \draw [domain=-4.3:7.06] plot(\x,{(--11.3184--3.92*\x)/2.96});
               \draw [domain=-4.3:7.06] plot(\x,{(--13.3932-3.94*\x)/2.68});
               \draw [domain=-4.3:7.06] plot(\x,{(-2.5436--0.02*\x)/-5.64});
               \begin{scriptsize}
               \draw [fill=qqqqff] (-2.54,0.46) circle (2.5pt);
               \draw [fill=qqqqff] (0.42,4.38) circle (2.5pt);
               \draw [fill=qqqqff] (3.1,0.44) circle (2.5pt);
               \draw [fill=qqqqff] (-1.0841432360742709,2.3880265251989385) circle (2.5pt);
               \draw [fill=qqqqff] (1.7397223641328283,2.439811151237558) circle (2.5pt);
               \draw [fill=qqqqff] (0.3798923608928009,0.4496457717698837) circle (2.5pt);
               \end{scriptsize}
               \end{tikzpicture}
            \caption{$ $ $\Ha(6,9), \alphahat(I)=\frac 9 4$}\label{fig: 94 s=6}
            \end{figure}
    \item[e)] $\alphahat(I)=16/7$ and $Z$ is the configuration $\Ha(7,17)$ of seven points
          indicated in Figure \ref{fig: 167 s=7},
          \begin{figure}[H]
            \centering
               \begin{tikzpicture}[line cap=round,line join=round,x=1.0cm,y=1.0cm,scale=0.7]
               \clip(-2,-1) rectangle (9,4);
               \draw [domain=-4.3:12.94] plot(\x,{(--1.7632--3.34*\x)/2.84});
               \draw [domain=-4.3:12.94] plot(\x,{(--13.0772-3.42*\x)/2.02});
               \draw [domain=-4.3:12.94] plot(\x,{(--1.6192--0.08*\x)/-4.86});
               \draw [domain=-0.1:12.94] plot(\x,{(--5.465860227434063-1.0747454838418833*\x)/1.9371067446167143});
               \begin{scriptsize}
               \draw [fill=qqqqff] (-0.8,-0.32) circle (2.5pt);
               \draw [fill=qqqqff] (2.04,3.02) circle (2.5pt);
               \draw [fill=qqqqff] (4.06,-0.4) circle (2.5pt);
               \draw [fill=qqqqff] (1.2715043805797765,2.1162058560339627) circle (2.5pt);
               \draw [fill=qqqqff] (3.208611125196491,1.0414603721920794) circle (2.5pt);
               \draw [fill=qqqqff] (8.062865995090155,-0.465890798273089) circle (2.5pt);
               \draw [fill=qqqqff] (5.860086440609344,-0.42963105252031847) circle (2.5pt);
               \end{scriptsize}
               \end{tikzpicture}
            \caption{$ $ $\Ha(7,17),\; \alphahat(I)=\frac{16}{7}$}\label{fig: 167 s=7}
            \end{figure}
    \item[f)] $\alphahat(I)=7/3$ and
             \begin{itemize}
               \item  $Z$ is one of the configurations $\Ha(6,2)$ or $\Ha(6,5)$ of six points indicated in Figure \ref{fig: 73 s=6}, or
                   \begin{figure}[H]
                    \centering
                       \begin{minipage}{0.4\textwidth}
                       \centering
                         \begin{tikzpicture}[line cap=round,line join=round,x=1.0cm,y=1.0cm,scale=0.55]
                         \clip(-4,-1.75) rectangle (7.5,4);
                         \draw [domain=-2.7:18.7] plot(\x,{(-1.854-0.8*\x)/-0.54});
                         \draw [domain=-1.7:2.6] plot(\x,{(--3.4696-1.64*\x)/1.4});
                         \draw [domain=-2.3:5.2] plot(\x,{(-3.318--0.84*\x)/-1.94});
                         \draw [domain=-4.3:18.7] plot(\x,{(-0.248-0.12*\x)/4.16});
                         \draw (0,-0.56) node[anchor=north west] {$\Ha(6,2)$};
                         \begin{scriptsize}
                         \draw [fill=qqqqff] (-0.36,2.9) circle (2.5pt);
                         \draw [fill=qqqqff] (-0.9,2.1) circle (2.5pt);
                         \draw [fill=qqqqff] (1.04,1.26) circle (2.5pt);
                         \draw [fill=qqqqff] (-0.68,-0.04) circle (2.5pt);
                         \draw [fill=qqqqff] (3.48,-0.16) circle (2.5pt);
                         \draw [fill=qqqqff] (5.79922401847575,-0.22690069284064662) circle (2.5pt);
                         \end{scriptsize}
                         \end{tikzpicture}
                       \end{minipage}
                       \quad
                       \begin{minipage}{0.4\textwidth}
                       \centering
                         \begin{tikzpicture}[line cap=round,line join=round,x=1.0cm,y=1.0cm,scale=0.55]
                         \clip(-2,-2) rectangle (5,4);
                         \draw (-0.62,-1.2) -- (-0.62,6.3);
                         \draw [domain=-4.3:7.06] plot(\x,{(-0.5712-0.*\x)/2.04});
                         \draw (0,-0.76) node[anchor=north west] {$\Ha(6,5)$};
                         \begin{scriptsize}
                         \draw [fill=qqqqff] (-0.62,-0.28) circle (2.5pt);
                         \draw [fill=qqqqff] (1.42,-0.28) circle (2.5pt);
                         \draw [fill=qqqqff] (-0.62,1.48) circle (2.5pt);
                         \draw [fill=qqqqff] (-0.62,2.9053082046195264) circle (2.5pt);
                         \draw [fill=qqqqff] (3.127205117952819,-0.28) circle (2.5pt);
                         \draw [fill=qqqqff] (0.86,3.18) circle (2.5pt);
                         \end{scriptsize}
                         \end{tikzpicture}
                       \end{minipage}
                    \vspace{0.3cm}
                    \caption{$ $ $\alphahat(I)=\frac 7 3$}\label{fig: 73 s=6}
                     \end{figure}
               \item $Z$ is one of the configurations  $\Ha(7,13)$, $\Ha(7,21)$ or $\Ha(7,29)$ of seven points indicated in Figure \ref{fig: 73 s=7}, or
                   \begin{figure}[H]
                    \centering
                       \begin{minipage}{0.4\textwidth}
                       \centering
                       \begin{tikzpicture}[line cap=round,line join=round,x=1.0cm,y=1.0cm,scale=0.7]
                       \clip(-2,-2) rectangle (5,4);
                       \draw (-0.86,-0.5) -- (-0.86,6.3);
                       \draw [domain=-4.3:7.06] plot(\x,{(--0.4144-0.02*\x)/1.66});
                       \draw [domain=-4.3:3.06]
                       plot(\x,{(--6.660710362844703-2.7168708272859217*\x)/3.060278664731495});
                       \draw (1.2,-1.1) node[anchor=north west] {$\Ha(7,13)$};
                       \begin{scriptsize}
                       \draw [fill=qqqqff] (-0.86,0.26) circle (2pt);
                       \draw [fill=qqqqff] (-0.86,1.62) circle (2pt);
                       \draw [fill=qqqqff] (0.8,0.24) circle (2pt);
                       \draw [fill=qqqqff] (-0.86,2.94) circle (2pt);
                       \draw [fill=qqqqff] (2.200278664731495,0.22312917271407834) circle (2pt);
                       \draw [fill=qqqqff] (3.9400261248185777,0.20216835994194482) circle (2pt);
                       \draw [fill=qqqqff] (1.2986731733973058,1.0235613965362738) circle (2pt);
                       \end{scriptsize}
                       \end{tikzpicture}
                       \end{minipage}
                       \quad
                       \begin{minipage}{0.4\textwidth}
                       \centering
                       \begin{tikzpicture}[line cap=round,line join=round,x=1.0cm,y=1.0cm,scale=0.55]
                       \clip(-3.72,-2) rectangle (4.88,5.78);
                       \draw [domain=-3.4:3.88] plot(\x,{(-5.0432-1.72*\x)/-0.98});
                       \draw [domain=-3.72:3.88] plot(\x,{(-6.7856--0.02*\x)/-2.24});
                       \draw [domain=-3.72:3.88] plot(\x,{(--5.6148-1.74*\x)/1.26});
                       \draw [domain=-0.85:1.88] plot(\x,{(-1.3274103577148955-1.726972770955686*\x)/-0.19904965296316066});
                       \draw [domain=-3.72:6.88] plot(\x,{(-0.11768838857619188--0.058724202626641464*\x)/6.177510944340213});
                       \draw (0.5,-0.76) node[anchor=north west] {$\Ha(7,21)$};
                       \begin{scriptsize}
                       \draw [fill=qqqqff] (-0.22,4.76) circle (2.5pt);
                       \draw [fill=qqqqff] (-1.2,3.04) circle (2.5pt);
                       \draw [fill=qqqqff] (1.04,3.02) circle (2.5pt);
                       \draw [fill=qqqqff] (-0.41904965296316066,3.033027229044314) circle (2.5pt);
                       \draw [fill=qqqqff] (-2.9589743589743596,-0.047179487179486834) circle (2.5pt);
                       \draw [fill=qqqqff] (3.218536585365854,0.01154471544715463) circle (2.5pt);
                       \draw [fill=qqqqff] (-0.7716755312971246,-0.026386747074474845) circle (2.5pt);
                       \end{scriptsize}
                       \end{tikzpicture}
                       \end{minipage}
                       \quad
                       \begin{minipage}{0.4\textwidth}
                       \centering
                       \begin{tikzpicture}[line cap=round,line join=round,x=1.0cm,y=1.0cm,scale=0.7]
                       \clip(-5.5,-1.34) rectangle (4.92,4.02);
                       \draw [rotate around={-177.3850944090728:(0.7738129265320349,0.9606544788109032)}] (0.7738129265320349,0.9606544788109032) ellipse (2.9242523364015285cm and 1.2553398245502152cm);
                       \draw [domain=-3.52:4.92] plot(\x,{(-3.576-1.92*\x)/-1.56});
                       \draw [domain=-3.52:4.92] plot(\x,{(--5.9856-2.06*\x)/1.26});
                       \draw [domain=-3.52:4.92] plot(\x,{(--2.266169578184917-3.6480528334043463*\x)/-0.2577929271410313});
                       \draw (-5,1.4) node[anchor=north west] {$\Ha(7,29)$};
                       \begin{scriptsize}
                       \draw [fill=qqqqff] (-0.14,2.12) circle (2pt);
                       \draw [fill=qqqqff] (-1.7,0.2) circle (2pt);
                       \draw [fill=qqqqff] (0.6,-0.3) circle (2pt);
                       \draw [fill=qqqqff] (2.82,0.14) circle (2pt);
                   \draw [fill=qqqqff] (1.56,2.2) circle (2pt);
                   \draw [fill=qqqqff] (0.8577929271410313,3.3480528334043465) circle (2pt);
                   \draw [fill=qqqqff] (0.7778810192347608,2.217211637358913) circle (2pt);
                   \end{scriptsize}
                   \end{tikzpicture}
                   \end{minipage}
                \caption{$ $ $\alphahat(I)=\frac 7 3$}\label{fig: 73 s=7}
                \end{figure}
               \item $Z$ is one of the configurations  $\Ha(8,119)$ or $\Ha(8,137)$ of eight points indicated in Figure \ref{fig: 73 s=8},
                   \begin{figure}[H]
                    \centering
                       \begin{minipage}{0.4\textwidth}
                       \centering
                       \begin{tikzpicture}[line cap=round,line join=round,x=1.0cm,y=1.0cm,scale=0.55]
                       \clip(-2.5,-4.5) rectangle (3.6,3.5);
                       \draw [domain=-1.9042372989941696:3.5804439496693936] plot(\x,{(--3.3464--0.76*\x)/2.02});
                       \draw [domain=0:2] plot(\x,{(--4.5344-3.68*\x)/1.02});
                       \draw [domain=-1.9042372989941696:3.5804439496693936] plot(\x,{(--2.7777584551663606--0.19049181227199186*\x)/2.177852812902872});
                       \draw [domain=-1.9042372989941696:3.5804439496693936] plot(\x,{(--1.9615243411870944-0.6269862901239636*\x)/2.4044364173713157});
                       \draw [domain=1.5:3.5804439496693936] plot(\x,{(-4.238895076497239--2.0119679836823217*\x)/0.4422216730475259});
                       \draw (0,-3.4) node[anchor=north west] {$\Ha(8,119)$};
                       \begin{scriptsize}
                       \draw [fill=qqqqff] (-1.32,1.16) circle (2.5pt);
                       \draw [fill=qqqqff] (1.72,-1.76) circle (2.5pt);
                       \draw [fill=qqqqff] (0.7,1.92) circle (2.5pt);
                       \draw [fill=qqqqff] (0.8578528129028717,1.3504918122719918) circle (2.5pt);
                       \draw [fill=qqqqff] (1.084436417371316,0.5330137098760364) circle (2.5pt);
                       \draw [fill=qqqqff] (2.162221673047526,0.25196798368232187) circle (2.5pt);
                       \draw [fill=qqqqff] (2.4339733009300373,1.4883514721831328) circle (2.5pt);
                       \draw [fill=qqqqff] (2.693719071039209,2.6701121257375244) circle (2.5pt);
                       \end{scriptsize}
                       \end{tikzpicture}
                       \end{minipage}
                       \quad
                       \begin{minipage}{0.4\textwidth}
                       \centering
                       \begin{tikzpicture}[line cap=round,line join=round,x=1.0cm,y=1.0cm,scale=0.55]
                       \clip(-1.5,-5) rectangle (13,3.7);
                       \draw [domain=-2.:7] plot(\x,{(--1.7632--3.34*\x)/2.84});
                       \draw [domain=-2.:4.5] plot(\x,{(--13.0772-3.42*\x)/2.02});
                       \draw [domain=-2.:10] plot(\x,{(--1.6192--0.08*\x)/-4.86});
                       \draw [domain=0.:7] plot(\x,{(--5.465860227434063-1.0747454838418833*\x)/1.9371067446167143});
                       \draw (0,-3.5) node[anchor=north west] {$\Ha(8,137)$};
                       \begin{scriptsize}
                       \draw [fill=qqqqff] (-0.8,-0.32) circle (2.5pt);
                       \draw [fill=qqqqff] (2.04,3.02) circle (2.5pt);
                       \draw [fill=qqqqff] (4.06,-0.4) circle (2.5pt);
                       \draw [fill=qqqqff] (1.2715043805797765,2.1162058560339627) circle (2.5pt);
                       \draw [fill=qqqqff] (3.208611125196491,1.0414603721920794) circle (2.5pt);
                       \draw [fill=qqqqff] (8.062865995090155,-0.465890798273089) circle (2.5pt);
                       \draw [fill=qqqqff] (5.860086440609344,-0.42963105252031847) circle (2.5pt);
                       \draw [fill=qqqqff] (9.660457800728011,-0.49218860577330065) circle (2.5pt);
                   \end{scriptsize}
                   \end{tikzpicture}
                   \end{minipage}
                \caption{$ $ $\alphahat(I)=\frac 7 3$}\label{fig: 73 s=8}
                \end{figure}
             \end{itemize}
    \item[g)] $\alphahat(I)=12/5$ and $Z$ is a set of six points not contained in a conic and no three points are collinear,
    \item[h)] $\alphahat(I)=17/7$ and
            \begin{itemize}
               \item $Z$ is the configuration $\Ha(8,136)$ of eight points indicated in Figure \ref{fig: h(8,136)}, or
                   \begin{figure}[H]
                    \centering
                    \begin{tikzpicture}[line cap=round,line join=round,>=triangle 45,x=1.0cm,y=1.0cm,scale=0.7]
                    \clip(-1.24,-0.96) rectangle (10.14,3.46);
                    \draw [domain=-1.24:10.14] plot(\x,{(--1.7632--3.34*\x)/2.84});
                    \draw [domain=-1.24:10.14] plot(\x,{(--13.0772-3.42*\x)/2.02});
                    \draw [domain=-1.24:10.14] plot(\x,{(--1.6192--0.08*\x)/-4.86});
                    \draw [domain=0.24:10.14] plot(\x,{(--5.465860227434063-1.0747454838418833*\x)/1.9371067446167143});
                    \begin{scriptsize}
                    \draw [fill=black] (-0.8,-0.32) circle (2.5pt);
                    \draw [fill=black] (2.04,3.02) circle (2.5pt);
                    \draw [fill=black] (4.06,-0.4) circle (2.5pt);
                    \draw [fill=black] (1.2715043805797765,2.1162058560339627) circle (2.5pt);
                    \draw [fill=black] (3.208611125196491,1.0414603721920794) circle (2.5pt);
                    \draw [fill=black] (7.300767967493439,-0.4533459747735546) circle (2.5pt);
                    \draw [fill=black] (8.639417590789808,-0.4753813595191738) circle (2.5pt);
                    \draw [fill=black] (9.840079911961396,-0.49514534834504365) circle (2.5pt);
                    \end{scriptsize}
                    \end{tikzpicture}
                    \caption{$ $ $\Ha(8,136),\;\;\alphahat(I)=\frac{17}{7}$}\label{fig: h(8,136)}
                    \end{figure}
               \item $Z$ is a configuration of nine points indicated in Figure \ref{fig: 9points}.
                   \begin{figure}[H]
                    \centering
                    \begin{tikzpicture}[line cap=round,line join=round,>=triangle 45,x=1.0cm,y=1.0cm,scale=0.7]
                    \clip(-1.16,-0.88) rectangle (10.2,3.4);
                    \draw [domain=-1.16:10.2] plot(\x,{(--1.9232--3.34*\x)/2.34});
                    \draw [domain=-1.16:10.2] plot(\x,{(--12.8772-3.42*\x)/2.52});
                    \draw [domain=-1.16:10.2] plot(\x,{(--1.6192--0.08*\x)/-4.86});
                    \draw [domain=0.16:10.2] plot(\x,{(--5.399713519566572-1.0747454838418833*\x)/2.0910683995749824});
                    \begin{scriptsize}
                    \draw [fill=black] (-0.8,-0.32) circle (2.5pt);
                    \draw [fill=black] (1.54,3.02) circle (2.5pt);
                    \draw [fill=black] (4.06,-0.4) circle (2.5pt);
                    \draw [fill=black] (0.9068029051255904,2.1162058560339627) circle (2.5pt);
                    \draw [fill=black] (2.9978713047005727,1.0414603721920794) circle (2.5pt);
                    \draw [fill=black] (7.300767967493439,-0.4533459747735546) circle (2.5pt);
                    \draw [fill=black] (5.860086440609344,-0.42963105252031847) circle (2.5pt);
                    \draw [fill=black] (8.639417590789808,-0.4753813595191738) circle (2.5pt);
                \draw [fill=black] (9.840079911961396,-0.49514534834504365) circle (2.5pt);
                \end{scriptsize}
                \end{tikzpicture}
                \caption{$ $ $\alphahat(I)=\frac{17}{7}$}\label{fig: 9points}
                \end{figure}
             \end{itemize}
\end{itemize}
\end{theorem}

\proof
   It suffices to show that the conditions
   $$\alphahat(I)< \frac52\;\; \mbox{ and }\;\; s\geq 10,$$
   imply $\alphahat(I)\leq 2$, and
   $$\frac73\leq\alphahat(I)< \frac52\;\; \mbox{ and }\;\; s=9,$$
   imply $\alphahat(I)= \frac{17}{7}$ and the configuration is exactly that
   indicated in Figure \ref{fig: 9points}. Indeed, the rest follows from looking up
   the list in \cite{Har}.

   Let $Z=\{P_1,\ldots, P_8, P_9, \ldots, P_s\}$ be a configuration of $s\geq 9$ points. Assume that
   $Z$ is not contained in a conic (as otherwise obviously $\alphahat(I)\leq 2$).
   Since $Z$ is not contained in a conic, there is also a subset $W\subseteq Z$ of $8$ points
   not contained in a conic. Renumbering the points if necessary, we may assume that $W=\{P_1,\ldots, P_8\}$.
   Then $\alphahat(I_W)=\frac 7 3$ or $\alphahat(I_W)=\frac{17}{7}$ and $W$ is
   one of the configurations $\Ha(8,119)$, $\Ha(8,136)$ or $\Ha(8,137)$.

\textbf{Case $\Ha(8,119)$}.
   In this case $W$ is indicated in Figure \ref{fig:case a}.
\begin{figure}[H]
\centering
\begin{tikzpicture}[line cap=round,line join=round,>=triangle 45,x=1.0cm,y=1.0cm]
\clip(-1.611365970182038,-1.995519983892074) rectangle (3.0213259583007774,3.1430406034480565);
\draw [domain=-1.611365970182038:3.0213259583007774] plot(\x,{(--3.3464--0.76*\x)/2.02});
\draw [domain=-1.611365970182038:3.0213259583007774] plot(\x,{(--4.5344-3.68*\x)/1.02});
\draw [domain=-1.611365970182038:3.0213259583007774] plot(\x,{(--2.7777584551663606--0.19049181227199186*\x)/2.177852812902872});
\draw [domain=-1.611365970182038:3.0213259583007774] plot(\x,{(--1.9615243411870944-0.6269862901239636*\x)/2.4044364173713157});
\draw [domain=-1.611365970182038:3.0213259583007774] plot(\x,{(-4.238895076497239--2.0119679836823217*\x)/0.4422216730475259});
\begin{scriptsize}
\draw [fill=black] (-1.32,1.16) circle (2pt);
\draw[color=black] (-1.251932975730785,1.5322482949813314) node {$P_2$};
\draw [color=black] (1.72,-1.76) circle (2pt);
\draw[color=black] (1.9962763074583156,-1.5029636581625798) node {$P_1$};
\draw [fill=black] (0.7,1.92) circle (2pt);
\draw [fill=black] (0.8578528129028717,1.3504918122719918) circle (2pt);
\draw [fill=black] (1.084436417371316,0.5330137098760364) circle (2pt);
\draw [fill=black] (2.162221673047526,0.25196798368232187) circle (2pt);
\draw [fill=black] (2.4339733009300373,1.4883514721831328) circle (2pt);
\draw [fill=black] (2.693719071039209,2.6701121257375244) circle (2pt);
\end{scriptsize}
\end{tikzpicture}
\caption{$ $ Case $\Ha(8,119)$}\label{fig:case a}
\end{figure}
   The idea now is to replace $W$ by $W'=\{P_2,\ldots, P_9\}$.
   Since $U=\{P_2,\ldots, P_8\}$ is not contained in a conic,
   neither is $W'$, so that it must be $\alphahat(I_{W'})=\frac 7 3$ and $W'$ is again one of the configurations  $\Ha(8,119)$, $\Ha(8,136)$ or $\Ha(8,137)$.
   Since $U$ has no $4$ collinear points, $W'=U\cup\{P_9\}$ cannot have $5$ collinear points,
   so that $W'$ must be again of the type $\Ha(8,119)$. Since the triple point of the configuration is fixed in $P_2$,
   it is elementary to check that the only way to complete $U$ to the
   $\Ha(8,119)$ configuration is by adding the point $P_1$, a contradiction.

\textbf{Case $\Ha(8,137)$}.
   In this case $W$ is indicated in Figure \ref{fig:case b}.

\begin{figure}[H]
\centering
\begin{tikzpicture}[line cap=round,line join=round,>=triangle 45,x=1.0cm,y=1.0cm]
\clip(-1.16,-1.02) rectangle (10.16,3.74);
\draw [domain=-5:10.16] plot(\x,{(--1.7632--3.34*\x)/2.84});
\draw [domain=-5:10.16] plot(\x,{(--13.0772-3.42*\x)/2.02});
\draw [domain=-5:10.16] plot(\x,{(--1.6192--0.08*\x)/-4.86});
\draw [domain=0.5:10.16] plot(\x,{(--5.465860227434063-1.0747454838418833*\x)/1.9371067446167143});
\begin{scriptsize}
\draw [fill=black] (-0.8,-0.32) circle (2.5pt);
\draw[color=black] (-0.92,0.22) node {$P_3$};
\draw [fill=black] (2.04,3.02) circle (2.5pt);
\draw[color=black] (2.04,3.62) node {$P_2$};
\draw [fill=black] (4.06,-0.4) circle (2.5pt);
\draw[color=black] (4.24,0.02) node {$P_4$};
\draw [fill=black] (1.2715043805797765,2.1162058560339627) circle (2.5pt);
\draw [fill=black] (3.208611125196491,1.0414603721920794) circle (2.5pt);
\draw[color=black] (3.38,1.46) node {$P_5$};
\draw [fill=black] (8.062865995090155,-0.465890798273089) circle (2.5pt);
\draw[color=black] (8.04,0.1) node {$P_6$};
\draw [color=black] (5.860086440609344,-0.42963105252031847) circle (2.5pt);
\draw[color=black] (6.04,0.) node {$P_1$};
\draw [fill=black] (9.660457800728011,-0.49218860577330065) circle (2.5pt);
\draw[color=black] (9.64,0.08) node {$P_7$};
\end{scriptsize}
\end{tikzpicture}
\caption{$ $ Case $\Ha(8,137)$}\label{fig:case b}
\end{figure}
   Now we replace $W$ by $W'=\{P_2,\ldots, P_9\}$.
   Since $U=\{P_2,\ldots,P_8\}$ is not contained in a conic, neither is $W'$.
   So again $W'$ must be one of the configurations $\Ha(8,119)$, $\Ha(8,136)$ or $\Ha(8,137)$.
   Assume that $W'$ is of type $\Ha(8,119)$. Since $U$ contains only one set of four collinear points,
   the point $P_9$ must lie either on the line $P_2 P_3$ or on the line $P_2 P_4$.
   Since these cases are equivalent we may assume, without loss of generality,
   that $P_9$ lies on the line $P_2 P_3$. But $P_9$ must also lie on two lines, each containing three points.
   There is only one such possibility: the point $P_9$ has to lie on both lines, $P_5 P_6$ and $P_5 P_7$, a contradiction.

   Now we consider the case that $W'$ is of type $\Ha(8,137)$.
   Then $P_9$ must be the intersection point of the line $P_3 P_4$ and a line passing through two other points.
   That must be again the point $P_1$. A contradiction.

   Hence $W'$ is of the type $\Ha(8,136)$. Then the point $P_9$ must be on the line $P_6 P_7$.
   Thus the set $Z'=W' \cup \{P_1\}$ is indicated in Figure \ref{fig: 9points}. We will
   complete the argument in the Case of $9$ points below.

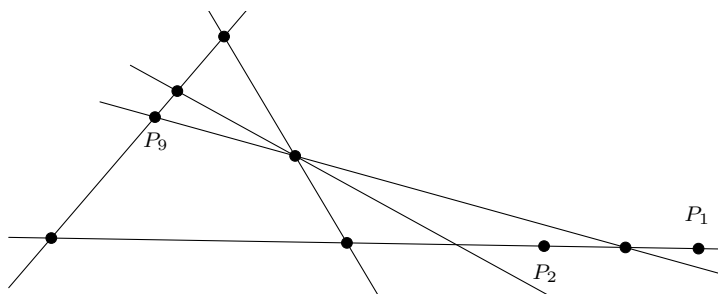
\begin{figure}[H]
\centering
\begin{tikzpicture}[line cap=round,line join=round,>=triangle 45,x=1.0cm,y=1.0cm,scale=0.8]
\clip(-1.5,-1.24) rectangle (10.16,3.44);
\draw [domain=-1.5:10.16] plot(\x,{(--1.7632--3.34*\x)/2.84});
\draw [domain=-1.5:10.16] plot(\x,{(--13.0772-3.42*\x)/2.02});
\draw [domain=-1.5:10.16] plot(\x,{(--1.6192--0.08*\x)/-4.86});
\draw [domain=0.5:10.16] plot(\x,{(--5.465860227434063-1.0747454838418833*\x)/1.9371067446167143});
\draw [domain=0:10.16] plot(\x,{(-10.522924978491007--1.5168417317112532*\x)/-5.430806465593317});
\begin{scriptsize}
\draw [fill=black] (-0.8,-0.32) circle (2.5pt);
\draw [fill=black] (2.04,3.02) circle (2.5pt);
\draw [fill=black] (4.06,-0.4) circle (2.5pt);
\draw [fill=black] (1.2715043805797765,2.1162058560339627) circle (2.5pt);
\draw [fill=black] (3.208611125196491,1.0414603721920794) circle (2.5pt);
\draw [fill=black] (7.300767967493439,-0.4533459747735546) circle (2.5pt);
\draw[color=black] (7.32,-0.9) node {$P_2$};
\draw [fill=black] (8.639417590789808,-0.4753813595191738) circle (2.5pt);
\draw [fill=black] (9.840079911961396,-0.49514534834504365) circle (2.5pt);
\draw[color=black] (9.82,0.06) node {$P_1$};
\draw [fill=black] (0.9047872142238311,1.6849258082773229) circle (2.5pt);
\draw[color=black] (0.92,1.28) node {$P_9$};
\end{scriptsize}
\end{tikzpicture}
\caption{$ $ The set $Z'$ }\label{fig: case d}
\end{figure}

\textbf{Case $\Ha(8,136)$}. In this situation $W$ is presented in Figure \ref{fig: case c}.
\begin{figure}[H]
\centering
\begin{tikzpicture}[line cap=round,line join=round,>=triangle 45,x=1.0cm,y=1.0cm,scale=0.7]
\clip(-1.12,-0.86) rectangle (10.28,3.36);
\draw [domain=-1.12:10.28] plot(\x,{(--1.7632--3.34*\x)/2.84});
\draw [domain=-1.12:10.28] plot(\x,{(--13.0772-3.42*\x)/2.02});
\draw [domain=-1.12:10.28] plot(\x,{(--1.6192--0.08*\x)/-4.86});
\draw [domain=0.5:10.28] plot(\x,{(--5.465860227434063-1.0747454838418833*\x)/1.9371067446167143});
\begin{scriptsize}
\draw [fill=black] (-0.8,-0.32) circle (2.5pt);
\draw [fill=black] (2.04,3.02) circle (2.5pt);
\draw [fill=black] (4.06,-0.4) circle (2.5pt);
\draw [fill=black] (1.2715043805797765,2.1162058560339627) circle (2.5pt);
\draw [fill=black] (3.208611125196491,1.0414603721920794) circle (2.5pt);
\draw [fill=black] (7.300767967493439,-0.4533459747735546) circle (2.5pt);
\draw [fill=black] (8.639417590789808,-0.4753813595191738) circle (2.5pt);
\draw [color=black] (9.840079911961396,-0.49514534834504365) circle (2.5pt);
\draw[color=black] (9.82,0.06) node {$P_1$};
\end{scriptsize}
\end{tikzpicture}
\caption{$ $ Case $\Ha(8,136)$}\label{fig: case c}
\end{figure}
   We proceed analogously and show that there are only two possibilities.
   The set $Z'=W'\cup \{P_9\}$ is either as in Figure \ref{fig: 9points},
   or as in Figure \ref{fig: case d}.

   If $Z'$ is the configuration indicated in Figure \ref{fig: case d}, then $Z''=\{P_1,P_3,\dots,P_9\}$ 
   is the configuration $\Ha(8,118)$ 
   and its Waldschmidt constant is $\alphahat(I_{Z''})=\frac52>\alphahat(I_Z)$, which is a contradiction.

   Now we turn to the final case.

\textbf{Case of $9$ points}.
   From what we have done so far, it follows that the configuration presented in Figure \ref{fig: 9points}
   is the only configuration of $9$ points with $2<\alphahat(I)<\frac52$.
   It is not hard to check that in this case $\alphahat(I)=\frac{17}{7}$. We verified hand calculations with Singular, \cite{DGPS}.

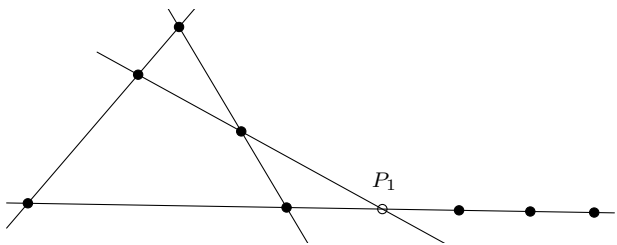
\begin{figure}[H]
\centering
\begin{tikzpicture}[line cap=round,line join=round,>=triangle 45,x=1.0cm,y=1.0cm,scale=0.7]
\clip(-1.2,-1.06) rectangle (10.22,3.36);
\draw [domain=-1.2:10.22] plot(\x,{(--1.7632--3.34*\x)/2.84});
\draw [domain=-1.2:10.22] plot(\x,{(--13.0772-3.42*\x)/2.02});
\draw [domain=-1.2:10.22] plot(\x,{(--1.6192--0.08*\x)/-4.86});
\draw [domain=0.5:10.22] plot(\x,{(--5.465860227434063-1.0747454838418833*\x)/1.9371067446167143});
\begin{scriptsize}
\draw [fill=black] (-0.8,-0.32) circle (2.5pt);
\draw [fill=black] (2.04,3.02) circle (2.5pt);
\draw [fill=black] (4.06,-0.4) circle (2.5pt);
\draw [fill=black] (1.2715043805797765,2.1162058560339627) circle (2.5pt);
\draw [fill=black] (3.208611125196491,1.0414603721920794) circle (2.5pt);
\draw [fill=black] (7.300767967493439,-0.4533459747735546) circle (2.5pt);
\draw [color=black] (5.860086440609344,-0.42963105252031847) circle (2.5pt);
\draw[color=black] (5.9,0.1) node {$P_1$};
\draw [fill=black] (8.639417590789808,-0.4753813595191738) circle (2.5pt);
\draw [fill=black] (9.840079911961396,-0.49514534834504365) circle (2.5pt);
\end{scriptsize}
\end{tikzpicture}
\caption{$ $ Case of $9$ points }\label{fig: 10 points}
\end{figure}

   Finally, we claim that there is no configuration of $s\geq 10$ points with $2<\alphahat(I)<\frac52$.
   Indeed, assume that such a configuration of $10$ points exists. Removing the point $P_{10}$ we obtain a configuration
   indicated in Figure \ref{fig: 10 points}. Its Waldschmidt constant is equal $\frac{17}{7}$.
   Removing further the point $P_1$ we obtain $\Ha(8,136)$ with the Waldschmidt constant equal to $\frac{17}{7}$.
   Adding back the point $P_{10}$ we obtain a configuration of $9$ points with $2<\alphahat(I)<\frac52$.
   The argument above shows that it must be $P_{10}=P_1$, a contradiction.

\endproof

\section{Initial sequences with first differences equal $2$}
   Now we are in the position to show how Theorem \ref{thm: conf of points}
   applies towards the proof of Theorem B.
\proofof{Theorem B}
    Assume to the contrary, that there exists a planar point configuration $Z=\{P_1,\ldots, P_s\}$
    with radical ideal $I=I_Z$ such that
    $$\alpha (I^{(m+i)})-\alpha(I^{(m+(i-1))})=2$$
    for $i=1,2,3,4$ and such that $Z$ is not contained in a conic. Then the inequality (\ref{inq: wc})
    immediately implies that
    $$\alphahat(I)\leq \frac 7 3.$$
    By Theorem A it must be then $s\leq 8$. On the other hand it must be $s\geq 6$.
    Hence we must be in one of the cases listed in Theorem \ref{thm: conf of points}.
    It follows from the discussion in \cite[Remark 5.11]{GHM09} that the sequence
    $\alpha(I^{(m)})$ in all these cases is determined by a finite number of initial
    terms. With some more effort one can in fact verify that in each of these cases
    the sequence of first differences
    $$\alpha(I^{(m)})-\alpha(I^{(m-1)}),\;\;\mbox{ with } \alpha(I^{(0)})=0$$
    of the initial sequence of $I$ is periodic with the period $\pi=(a_1,\ldots,a_\ell)$ of length $\ell$
    of one of the following shapes:
    \begin{itemize}
           \item  $\Ha(6,10)$ : $\ell=2$, $\pi=(3,1)$,
           \item  $\Ha(6,2)$, $\Ha(6,5)$, $\Ha(7,13)$, $\Ha(7,21)$, $\Ha(7,29)$, $\Ha(8,119)$, $\Ha(8,137)$ : $\ell=3$, $\pi=(3,2,2)$,
           \item  $\Ha(6,9)$ : $\ell=4$, $\pi=(3,2,2,2),$
           \item  $\Ha(7,17)$ : $\ell=7$, $\pi=(3,2,2,3,2,2,2).$
    \end{itemize}
    We omit easy but tedious calculations. Since in none of the listed cases
    a subsequence of $4$ consecutive $2$'s appears, this contradicts the initial
    assumption and completes the proof of Theorem B.
\endproof

\begin{remark}
    Note that the configuration $\Ha(6,9)$ shows that three consecutive differences of $2$
    are possible without forcing $Z$ to be contained in a conic. This was already observed in
    the discussion in
    \cite[Lemma 4.12]{DST13}.
\end{remark}

The bound $\alphahat(I) < \frac52$ taken on in Theorem \ref{thm: conf of points}
might appear a bit incidental. However, already the list in \cite{Har} exhibits
a considerable number of cases with $\alphahat(I) = \frac52$. Moreover, a number
of experiments we have run for $s \geq9$ points show that it might be
significantly more difficult to chart all configuration of points with
Waldschmidt constants $\alphahat(I) \geq \frac52$ and the problem might
require some new methods. On the other hand, in the view of results obtained
so far it is tempting to conclude this note with the following question.

\noindent
\begin{problem} Is the set of Waldschmidt constants of all planar point
configurations well-ordered?
\end{problem}

\paragraph*{Acknowledgement.}
   These notes originated from discussions held during two workshops
   organized at the Banach Center in B\c edlewo and in Lanckorona. The final version of the paper
   was prepared during the Simons Semester on algebraic geometry in IMPAN
   and was partially supported by the grant 346300 for IMPAN from the Simons Foundation and the matching 2015-2019 Polish MNiSW fund.\\
   Additionally, the  research of Szpond was partially supported by National Science Centre, Poland, grant
   2014/15/B/ST1/02197.\\
   We would like to thank Marcin Dumnicki, Brian Harbourne and Tomasz Szemberg for helpful conversations.


\bigskip \small

\bigskip
   \L ucja Farnik,
   Jagiellonian University, Faculty of Mathematics and Computer Science, {\L}ojasiewicza 6, PL-30-348 Krak\'ow, Poland

\nopagebreak
   \textit{E-mail address:} \texttt{lucja.farnik@gmail.com}

\bigskip
   Janusz Gwo\'zdziewicz,
   Department of Mathematics, Pedagogical University of Cracow,
   Podchor\c a\.zych 2,
   PL-30-084 Krak\'ow, Poland

\nopagebreak
   \textit{E-mail address:} \texttt{gwozd63@gmail.com}

\bigskip
   Beata Hejmej,
   Department of Mathematics, Pedagogical University of Cracow,
   Podchor\c a\.zych 2,
   PL-30-084 Krak\'ow, Poland

\nopagebreak
   \textit{E-mail address:} \texttt{bhejmej1f@gmail.com}

\bigskip
   Magdalena Lampa-Baczy\'nska,
   Department of Mathematics, Pedagogical University of Cracow,
   Podchor\c a\.zych 2,
   PL-30-084 Krak\'ow, Poland

\nopagebreak
   \textit{E-mail address:} \texttt{lampa.baczynska@wp.pl}

\bigskip
   Grzegorz Malara,
   Department of Mathematics, Pedagogical University of Cracow,
   Podchor\c a\.zych 2,
   PL-30-084 Krak\'ow, Poland

\nopagebreak
   \textit{E-mail address:} \texttt{grzegorzmalara@gmail.com}

\bigskip
   Justyna Szpond,
   Department of Mathematics, Pedagogical University of Cracow,
   Podchor\c a\.zych 2,
   PL-30-084 Krak\'ow, Poland

\nopagebreak
   \textit{E-mail address:} \texttt{szpond@gmail.com}


\end{document}